\documentclass[12pt,leqno]{amsart}
\usepackage{amsmath,amssymb,amsfonts}
\usepackage{eucal}
\usepackage{graphicx}

\newcommand \comment[1]{}			
\newcommand \dateadded[1]{\comment{[Date added: #1.]}}
\newcommand \mylabel[1]{\label{#1}\comment{{\rm \{#1\} }}}
\newcommand \myref[1]{\ref{#1}\comment{{\{#1\}}}}

\newtheorem{lem}{Lemma}

\newtheorem{prop}[lem]{Proposition}
\newtheorem{thm}[lem]{Theorem}

\theoremstyle{remark}
\newtheorem{exam}{Example}

\newcommand \myexam[1]{\smallskip\begin{exam}[\emph{#1}]}
\newcommand \myprob[1]{\smallskip\begin{prob}[\emph{#1}]}

\renewcommand\phi{\varphi}                 
\renewcommand\epsilon{\varepsilon}
\newcommand\eset{\varnothing}
\newcommand\setm{\smallsetminus}

\newcommand \0{{\hat0}}
\newcommand\cA{{\mathcal A}}
\newcommand\cH{{\mathcal H}}
\newcommand\cL{{\mathcal L}}
\newcommand\cP{{\mathcal P}}
\newcommand\cS{{\mathcal S}}
\newcommand\bbP{{\mathbb P}}
\newcommand\bbR{{\mathbb R}}

\begin{document}


\begin{center}

\large
\textsc{On the division of space by topological hyperplanes}
\normalsize
\vskip20pt

{David Forge\\
Laboratoire de recherche en informatique UMR 8623\\
B\^at.\ 490, Universit\'e Paris-Sud\\
91405 Orsay Cedex, France\\
E-mail: {\tt forge@lri.fr}}\\[10pt]

and\\[10pt] 

{Thomas Zaslavsky\\
Department of Mathematical Sciences \\
State University of New York at Binghamton\\
Binghamton, NY 13902-6000, U.S.A.\\
E-mail: {\tt zaslav@math.binghamton.edu}}\\[20pt]

{Version of \today.}\\[20pt]

\end{center}

\small
{\sc Abstract.}
A \emph{topological hyperplane} is a subspace of $\bbR^n$ (or a homeomorph of it) that is topologically equivalent to an ordinary straight hyperplane.  
An \emph{arrangement} of topological hyperplanes in $\bbR^n$ is a finite set $\cH$ such that $k$ topological hyperplanes in $\cH$, if their intersection is nonempty, meet in a subspace that is a topological hyperplane in the intersection of any $k-1$ of them; but two topological hyperplanes that do intersect need not cross each other.  If every intersecting pair does cross, the arrangement is \emph{affine}.  
The number of regions formed by an arrangement of topological hyperplanes has the same formula as for arrangements of affine hyperplanes.
Hoping to explain this geometrically, we ask whether parts of the topological hyperplanes in any arrangement can be reassembled into an arrangement of affine topological hyperplanes with the same regions.  That is always possible if the dimension is two but not in higher dimensions.  
We also ask whether all affine topological hyperplane arrangements correspond to oriented matroids; they need not, but we can characterize those that do if the dimension is two.  In higher dimensions this problem is open.  Another open question is to characterize the intersection semilattices of topological hyperplane arrangements; a third is to prove that the regions of an arrangement of topological hyperplanes are necessarily cells.

\bigskip



 \normalsize

\date{Version of \today.}

\vskip 20pt



\section*{Dedication}

We dedicate this paper to Michel Las Vergnas on the occasion of his 65th birthday.  Not only is it related to his research into oriented matroids and counting via the Tutte polynomial, but the problem and its solution occurred to us during the meeting in his honor at the CIRM, Marseille-Luminy (whose generous hospitality we greatly enjoyed), in November of 2005.

\section{Introduction}

In a topological space $X$ that is homeomorphic to $\bbR^n$, a \emph{topological hyperplane}, or \emph{topoplane} for short, is a subspace $Y$ such that $(X,Y)$ is homeomorphic to $(\bbR^n,\bbR^{n-1})$.  
Consider a finite set $\cH$ of topoplanes in $X$.  Its \emph{intersection semilattice} is the class 
$$
\textstyle   \cL := \{ \bigcap\cS : \cS \subseteq \cH \text{ and } \bigcap\cS \neq \eset \},
$$ 
partially ordered (as is customary) by reverse inclusion; the members of $\cL$ are called the \emph{flats} of $\cH$.  
We study the combinatorial topology of an \emph{arrangement of topoplanes} in $X$, which is a finite set $\cH$ of topoplanes such that, for every topoplane $H\in\cH$ and flat $Y\in\cL$, either $Y\subseteq H$ or $H \cap Y$ is a topoplane in Y.  
We find that the simplest structure appears only in the planar case.  (There we call a topoplane a \emph{topological line}, abbreviated to \emph{topoline}.)

Zaslavsky showed in \cite[Theorem 3.2(A)]{CATD} that the number of \emph{regions} of a topoplane arrangement $\cH$---these are the components of the complement, $X \setm \bigcup\cH$---equals 
\begin{equation}\mylabel{E:count}
\sum_{Y \in \cL} |\mu(\0,Y)|, 
\end{equation}
where $\mu$ is the M\"obius function of $\cL$ and $\0$ is the zero element of $\cL$, that is, $X$.  The proof combined topology with combinatorics and assumed the side condition that every region is a topological cell.  Our work was inspired by the hope that, in a sense, Equation \eqref{E:count} would be no more general than the widely known formula for the number of regions of an arrangement of pseudospheres, or equivalently, topes of an oriented matroid.  We hoped, in particular, that the parts of the topoplanes of any  arrangement could be reorganized so that any two topoplanes that intersect actually cross, while not only the number but the actual regions remain exactly the same, and moreover that the reorganized arrangement is equivalent to an arrangement of pseudohyperplanes that represents an oriented matroid.  This hope, alas, failed, except in the plane.  
Even there, not every topoline arrangement represents an oriented matroid; but it is easy to characterize those that do (see Theorem \myref{T:2dim}).

The technical definition of crossing of topoplanes $H_1,H_2\in\cH$ is that 
\begin{equation}\mylabel{E:affdef}
(X,H_1,H_2,H_1\cap H_2) \cong (\bbR^n,\{x_1=0\},\{x_2=0\},\{x_1=x_2=0\}).
\end{equation}
We say $H_1$ and $H_2$ \emph{cross} and we call an arrangement \emph{affine} if every pair of topoplanes is disjoint or crossing.  
Our main theorem is that, for every arrangement $\cH$ of topolines, there is an affine topoline arrangement $\cA$ such that $\bigcup\cA = \bigcup\cH$.  

In the enumerative sense the least complicated topoplane arrangements $\cA$ are those that realize an oriented matroid.  Combinatorially, this means the regions are cells that correspond to the topes of an oriented matroid on the ground set $\cA$ (and this entails that the arrangement is affine) \cite{FL,OM}; thus the region-counting formula becomes the known formula for the number of topes (still assuming all the regions are cells).  Topologically, it means $\cA$ is isotopic to the affine part of a projective pseudohyperplane arrangement $\cP$ (which we will explain later).  In two dimensions, this is true given the obvious necessary condition, that the union $\bigcup\cA$ be connected, is sufficient; but in higher dimensions it is hopelessly far from the facts.  
Finding a necessary and sufficient condition for an affine topoplane arrangement $\cA$ whose union is connected to realize an oriented matroid is one open question.  
A second is whether the regions of a topoplane arrangement are necessarily open cells (as is known to be true for arrangements that realize an oriented matroid; see \cite{EM,Mandel} as described in \cite[p.\ 227]{OM}).  We expect that they must be, but we do not prove it.

One more question, that might turn out to be interesting, is to characterize the intersection semilattice.  We can prove each interval is a geometric lattice.  Though the intersection semilattice is not necessarily a geometric semilattice \cite{WW}, could it be true that every geometric semilattice is the intersection semilattice of an arrangement of topoplanes?

\section{Elementary properties}

We regard arrangements as topological objects, so we have to define homeomorphism.  We call two topoplane arrangements, $\cA$ in $X$ and $\cA'$ in $X'$, \emph{homeomorphic} if there is a homeomorphism $X \to X'$ that induces homeomorphisms of the topoplanes and consequently of all the flats and faces of the two arrangements.  
(An especially good kind of homeomorphism of arrangements in $X$ is \emph{isotopy}, which means a homotopy that at every stage (every value of $t\in[0,1]$) is a homeomorphism; isotopy is not central to this work but we mention it in connection with planar projectivization at Theorem \myref{T:2dim}.)

If $\cH$ is a topoplane arrangement, a flat $Y$ \emph{induces} the set
$$
\cH^Y := \{ Y\cap H : H \in \cH \text{ and } Y \not\subseteq H \text{ and } Y\cap H \neq \eset \}
$$
of topological subspaces of $Y$.  

\begin{prop}\mylabel{P:hereditary}
If $\cH$ is an arrangement of topoplanes and $Y$ is a flat, then the induced collection $\cH^Y$ is an arrangement of topoplanes.
\dateadded{23 Aug 2006}
\end{prop}

\begin{proof}
It is clear that $\cL(\cH^Y) = \{ Z \in \cL(\cH): Z \subseteq Y \}.$  This makes the lemma obvious from the definition.
\end{proof}

We often call an element of $\cH^Y$ a \emph{relative topoplane} in $Y$.

\begin{prop}\mylabel{P:geometric}
For an arrangement of topoplanes, each interval in $\cL$ is a geometric lattice with rank given by codimension.
\dateadded{24 Aug 2006}
\end{prop}

\begin{proof}
Consider a lower interval $[X,Y]$ in the partial ordering.  In this interval no two flats are disjoint.  Consequently, the function $r(Z) := \dim X - \dim Z$ is well defined and, since by definition $H \supseteq Z$ or $\dim(H\cap Z) = \dim Z - 1$ for any topoplane $H$ and flat $Z$ in the interval, $r$ satisfies the axioms of the rank function of a geometric lattice.
\end{proof}

In a topoplane arrangement there are only three possible relationships between two topoplanes $H_1$ and $H_2$.  They can be disjoint, they can cross as in \eqref{E:affdef}, or they can touch without crossing as in the next lemma.  In $\bbR^n$ let $G_+ := \{ x: x_1x_2=0 \text{ and } x_1, x_2 \geq 0\}$ and $G_- := \{ x: x_1x_2=0 \text{ and } x_1, x_2 \leq 0\}.$  Each of these sets is a topoplane that is the union of two perpendicular half-hyperplanes; their union is the union of the first two coordinate hyperplanes; and their intersection is the $x_1=x_2=0$ coordinate flat.

\begin{lem}\mylabel{L:noncrossing}
If topoplanes $H_1$ and $H_2$ in a topoplane arrangement $\cH$ meet but do not cross, then
\begin{equation}\mylabel{E:noncrossing}
(X,H_1,H_2,H_1\cap H_2) \cong (\bbR^n, G_+, G_-, \{x_1=x_2=0\}).
\end{equation}
\dateadded{24 Aug 2006}
\end{lem}

\begin{proof}
The intersection $Z := H_1 \cap H_2$ is a relative topoplane in each $H_i$ so it divides $H_i$ into two halves, the components of $H_i \setm Z$.  Call these halves $X_i^+$ and $X_i^-$.  Each half of $H_1$ lies in one of the (open) halfspaces formed by $H_2$.  Either both are in different halfspaces and we have Equation \eqref{E:affdef}, or they are in the same halfspace and we have Equation \eqref{E:noncrossing}.
\end{proof}

To clarify the idea of an affine topoplane arrangement we like to have a second characterization.  

\begin{prop}\mylabel{P:affdef}
A topoplane arrangement $\cH$ is affine if and only if, for each pair $H_1,H_2\in\cH$, each of the four regions into which they divide $X$ has boundary that intersects both $H_1 \setm H_2$ and $H_2 \setm H_1$.  
\dateadded{23 Aug 2006}
\end{prop}

\begin{proof}
This is obvious from Lemma \myref{L:noncrossing}.
\end{proof}

There is a more specific version of the characterization.

\begin{lem}\mylabel{L:crossing}
Topoplanes $H_1$ and $H_2$ cross if and only if they intersect each other and each of the regions they form has boundary that meets both $H_1 \setm H_2$ and $H_2 \setm H_1$.
\dateadded{23 Aug 2006}
\end{lem}

\begin{proof}
This also is obvious from Lemma \myref{L:noncrossing}.
\end{proof}

It will help us to have a general conception of crossing that we can apply to half topoplanes as well as whole ones.  Suppose $M$ is a manifold in $X$ and $H$ is a topoplane.  We say $H$ and $M$ \emph{cross} if they intersect and, at each intersection point, every neighborhood contains an open neighborhood $U$ such that $(U, H\cap U) \cong (\bbR^n, \{x_1=0\})$ and $M\cap U$ meets both components of $U\setm H$.  It is clear that this definition generalizes that given in the introduction, where $M$ is a topoplane.

\begin{lem}\mylabel{L:nbdcross}
Assume $\cH$ is a topoplane arrangement, $H \in \cH$, $Y \in \cL$, and $Z \in \cH^Y$ such that $Z \not\subseteq H$. 
Let $Z_+$ be either of the components of $Z\setm H$.  
Then $H \cap Z_+$ is a topoplane in $Z_+$ and $H$ and $Z$ cross if and only if $H$ and $Z_+$ cross.
\dateadded{23 Aug 2006}
\end{lem}

\begin{proof}
The first statement is obvious and the second is immediate from Lemma \myref{L:noncrossing}.
\end{proof}

\begin{lem}\mylabel{L:hereditarycrossing}
If in a topoplane arrangement $\cH$ two topoplanes, $H_1$ and $H_2$, cross, then $Y\cap H_1$ and $Y\cap H_2$ cross in $\cH^Y$ for each $Y \in \cL$ such that $Y \not\subseteq H_1, H_2$, both $Y\cap H_1$ and $Y\cap H_2$ are nonvoid, and $Y\cap H_1, Y\cap H_2$ are distinct.
\dateadded{23 Aug 2006}
\end{lem}

\begin{proof}
Suppose $Y$ has codimension $1$ and two relative topoplanes in $\cH^Y$ intersect.  
The relative topoplanes have the form $Y\cap H_1$ and $Y\cap H_2$ for $H_1, H_2 \in \cH$, and their intersection is $W := Y \cap Z$ where $Z := H_1 \cap H_2$.  
The set $Z_1 := Y\cap H_1$ cannot be in $H_2$, or else $Y\cap H_1 = Y\cap H_2$, contrary to the hypothesis that we have two different relative topoplanes; similarly $Z_2 := Y\cap H_2$ cannot be in $H_1$.  
Thus, $W$ has dimension $n-3$ by Proposition \myref{P:geometric}.  In $Y$ we have the relative topoplanes $Z_1$ and $Z_2$ whose intersection is $W$, a relative topoplane of both.  By Lemma \myref{L:noncrossing}, $Z_1$ and $Z_2$ form four regions in $Y$.  Each of these is the intersection with $Y$ of a different region of $\{H_1,H_2\}$ in $X$.

Let $R_+$ and $R_-$ be the regions of $\{H_1\}$ and let $S_+$ and $S_-$ be the regions of $\{H_2\}$.  Then $R_{ij} := R_i \cap S_j$ are the four regions of $\{H_1,H_2\}$.  
The intersections $Y \cap R_{ij}$ are the four regions of $\{Z_1,Z_2\}$ in $Y$.  What separates $Y \cap R_{++}$ from $Y \cap R_{+-}$ is $Y \cap H_2 = Z_2$, just as $H_2$ separates $R_{++}$ from $R_{+-}$ in $X$.  Similarly, $Z_1$ separates $Y \cap R_{++}$ from $Y \cap R_{-+}$.  This shows that $Z_1$ and $Z_2$ are both on the boundary of $Y \cap R_{++}$.  Similarly, both relative topoplanes are on the boundary of each $Y \cap R_{ij}$.  By Lemma \myref{L:crossing}, $Z_1$ and $Z_2$ cross in $Y$.

If $Y$ has codimension $d>1$, we apply induction on a maximal chain $Y \subset Y_1 \subset \cdots \subset Y_d = X$.
\end{proof}

\begin{prop}\mylabel{P:hereditaryaffine}
If $\cH$ is an affine arrangement of topoplanes and $Y$ is a flat, then so is the induced arrangement $\cH^Y$.
\dateadded{23 Aug 2006}
\end{prop}

\begin{proof}
We appeal to the previous lemma.
\end{proof}

\section{Reglueing}

The basic theorem is that, as concerns its combinatorics, a topoplane arrangement can be replaced by an affine topoplane arrangement.  
A \emph{face} of an arrangement is a region of the arrangement induced in a flat.  Thus, a $k$-dimensional face is a region of $\cH^t$ where $t$ is a $k$-dimensional flat of $\cH$.  A region of $\cH$ is a $d$-dimensional face where $d = \dim X$.  
The \emph{$k$-skeleton} of $\cH$ is the union of all $k$-dimensional flats.  Thus, writing $\cH^k$ for the $k$-skeleton, the $k$-faces are the components of $\cH^k \setm \cH^{k-1}$.  

\begin{thm}\mylabel{T:planar}
For any arrangement of topolines, there is an affine topoline arrangement which has the same faces.
\dateadded{23--25 Aug 2006}
\end{thm}

\begin{proof}
We apply the method of descent to the number of noncrossing intersecting pairs of topolines.  Suppose we have a noncrossing pair of topolines that intersect.  Their intersection $Z$ lies in $k\geq2$ topolines, call them $H^1, H^2, \ldots, H^k$.  $Z$ separates $H^i \setm Z$ into two halves, $H^i_+$ and $H^i_-$.  In cyclic order around $Z$, call these $2k$ halves $K^1_+, K^2_+, \ldots,K^k_+, K^1_-, K^2_-, \ldots, K^k_-$.  Let $K^i = K^i_+ \cup K^i_-$.  

It is clear that the new arrangement $\cH'$, which is $\cH$ with $H^1, \ldots, H^k$ replaced by $K^1, \ldots, K^k$, has the same skeleton in each dimension, hence it has the same faces.  
However, we have to check that $\cH'$ is an arrangement of topolines, and then that it has fewer noncrossing pairs of topolines than did $\cH$.

To show that $\cH'$ is an arrangement of topolines we consider the intersection of a topoline $H$ and a flat $Y$ of $\cH'$.  If $Y$ and $H$ are comparable or disjoint, the definition of a topoline arrangement is satisfied.  The only other case is that of two topolines.  If they both contain $Z$, they intersect in $Z$, which is a relative topoplane of both.  If neither contains $Z$, they are common topolines of $\cH$ and $\cH'$ so their intersection remains the same as in $\cH$.  Suppose the topolines are $H \not\supseteq Z$ and $K^1$ and suppose that $H \cap K^1$ consists of more than one point.  Then it consists of a point $W_+ \in K^1_+$ and a point $W_- \in K^1_-$.  
$K^1$ divides the plane into halves, $K^{1+}$ and $K^{1-}$, with $K^i_+$ in $K^{1+}$ for  $i = 2,\ldots,k$.  Also, $H$ divides the plane into two halves, $H^+$ and $H^-$; by choice of notation assume $O \in H^-$ and that the segment of $H$ from $W_+$ to $W_-$ lies in $K^{1+}$.  (All this is just to fix the notation.)  

Now, observe that $K^i_+$ is a topoline in $K^{1+}$ by Lemma \myref{L:nbdcross}.  It follows that $H$ intersects $H^i_+$.  Thus, $H$ intersects more than $k$ of the $2k$ half-topolines $H^i_\epsilon$, and consequently $H$ must intersect a topoline $H^i$ of $\cH$ more than once.  This is contrary to hypothesis, so it is impossible after all for $H \cap K^1$ to have more than one point.  The argument applies equally to each $K^i$, so we may conclude that $\cH'$ is a topoline arrangement.

Finally, we prove that the number of noncrossing pairs of topolines decreases from $\cH$ to $\cH'$.   A crossing pair from $\cH$, neither of them an $H^i$, remains crossing.  Amongst the $H^i$, the number of crossing pairs increases.  Suppose, then, that $H$ crosses exactly $k$ of the $H^i$, where $H \not\supseteq Z$.  Then $H$ crosses exactly $2k$ of the halves $K^i_+$ and $K^i_-$; hence by Lemma \myref{L:nbdcross} it crosses $k$ of the new topoplanes $K^i$.  Consequently, the number of crossing pairs increases.

Since there are fewer noncrossing topoline pairs in the new arrangement, by continuing the process we get an affine topoline arrangement.
\end{proof}

Reglueing can be impossible for a topoplane arrangement in three or more dimensions.  We give an example of this.

\myexam{Failure in three dimensions}\mylabel{X:3}
The example $\cH$ has five topoplanes in $\bbR^3$.  They are: 
\begin{align*}
H_1 &= \{ x: x_1 = 0 \} , \\
H_2 &= \{ x: x_2 = 0 \} , \\
H_3 &= \{ x: x_2 = |x_1| \} , \\
H_4 &= \{ x: x_3 = 0 \} , \\
H_5 &= \{ x: x_2+x_3 = 0 \} .
\end{align*}
Every pair crosses except $H_2$ and $H_3$.  
The common point of all topoplanes is $O$, the origin.  The $1$-dimensional flats are:
\begin{align*}
Z := H_1 \cap H_2 \cap H_3 &= \{ x: x_1 = x_2 = 0 \} , \\
H_1 \cap H_4 &= \{ x: x_1 = x_3 = 0 \} , \\
H_1 \cap H_5 &= \{ x: x_1  = 0, x_2+x_3 = 0 \} , \\
Y := H_2 \cap H_4 \cap H_5 &= \{ x: x_2 = x_3 = 0 \} , \\
H_3 \cap H_4 &= \{ x: x_2 = |x_1| , x_3 = 0 \} , \\
H_3 \cap H_5 &= \{ x: x_2 = |x_1| = -x_3 \} .
\end{align*}
The only two $1$-dimensional flats that lie in three topoplanes are $Z$ and $Y$.  This so limits the possibilities of recombining the faces of $\cH$ that it is impossible to get an affine arrangement $\cH'$.

To see why, note that $Y$ and $Z$ are relative topoplanes in a plane; therefore, in an affine recombination they have to cross.  This means, in effect, that they cannot be changed.  The plane $H_1$ that contains both has to remain a plane in $\cH'$.  Hence, the only potential changes in topoplanes are that $H_1$ and $H_3$ might be recombined and $H_4$ and $H_5$ might be recombined.  However, 
there is no way to recombine the halves of $H_1$ and $H_3$ so that two halves are on each side of $H_2$, which is a necessity if the recombined planes are to cross $H_2$.
\dateadded{25 Aug 2006}
\end{exam}

An intersection flat is \emph{simple} if its codimension equals the number of topoplanes that contain it; otherwise it is \emph{multiple}.  
It is no coincidence that our counterexample has multiple intersections.  
We call an arrangement \emph{simple} if every flat is simple.

\begin{thm}\mylabel{T:reglueing}
For a simple topoplane arrangement, there is an affine topoplane arrangement which has the same faces.
\dateadded{23--31 Aug 2006}
\end{thm}

\begin{proof}
The method of proof is similar to that of Theorem \myref{T:planar}, applying the method of descent to the number of noncrossing intersecting pairs of topoplanes.  

Suppose we have two noncrossing topoplanes, $H^1$ and $H^2$.  Their intersection $Z$ lies in no other topoplanes than these two.  $Z$ separates $H^i \setm Z$ into two halves.  In cyclic order around $Z$, call these four halves $H^1_+=K^1_+,\ H^2_+=K^2_+,\ H^2_-=K^1_-,\ H^1_-=K^2_-$, and let $K^i = K^i_+ \cup K^i_-$.  

The new arrangement $\cH'$, which is $\cH$ with $H^1, H^2$ replaced by $K^1, K^2$, has the same faces as $\cH$.  We need to prove that $\cH'$ is an arrangement of topoplanes and that it has fewer noncrossing pairs of topoplanes.

To show that $\cH'$ is an arrangement of topoplanes we consider the intersection of a topoplane $H$ and a flat $Y$ of $\cH'$.  There are four cases, depending mostly on whether either of them is a topoplane or flat in $\cH$.  

Before we can treat the cases we need to understand the flats of $\cH'$.  Those that are contained in $Z$, and those that are not contained in any $K^i$, are flats of $\cH$ because they are the intersection of topoplanes common to $\cH$ and $\cH'$.  Any other flat $V$ is the intersection of one $K^i$ with a flat $W$ not contained in either $K^1$ or $K^2$; so $W$ is a common flat of $\cH$ and $\cH'$.  Then 
\begin{equation}\mylabel{E:reflat}
V = V_+ \cup V_- \cup (W \cap Z), \text{ where } V_+ := W \cap K^i_+ \text{ and } V+- := W \cap K^i_-.
\end{equation}
Each $V_\epsilon$ is an intersection $W \cap H^j_\epsilon$.  Thus, it has codimension $1$ in $W$.  It follows that $V$ is a relative topoplane in $W$, assembled from the two half flats $V \cap H^1_\epsilon$ and $V \cap H^2_\epsilon$ as well as $V \cap Z$.  
 
Now we analyze the cases.  When $Y \in \cL$ (Cases 1 and 2), either $Y \subseteq Z$ or $Y \not\subseteq K^1, K^2$.  When $Y \notin \cL$ (Cases 3 and 4) we may assume $Y \subseteq K^2$ but $Y \not\subseteq K^1$.

\emph{Case 1.}  If $Y \in \cL$ and $H \neq K^1, K^2$, then $H\cap Y$ is empty or it is in $\cL$, hence is $Y$ or a relative topoplane of $Y$.  

\emph{Case 2.}  Suppose $Y \in \cL$ and $H=K^1$.  
If $Y \subseteq K^1$, then $Y \cap H = Y$. 
If $Y \not\subseteq K^1, K^2$, then $Y \cap H$ has the form of $V$ in \eqref{E:reflat} with $i=1$ and $W=Y$.  Thus, $Y \cap H$ is a relative topoplane in $Y$.

\emph{Case 3.}  Suppose $Y \notin \cL$ (so we assume $Y \subseteq K^2$ but $Y \not\subseteq K^1$) and $H=K^1$, then $Y$ has the form of $V$ in \eqref{E:reflat} with $i=2$.  Then $Y \cap H = Y \cap Z$, which is a relative topoplane in $Y$, as \eqref{E:reflat} shows.

\emph{Case 4.}  If $Y \not\in \cL$ and $H=K^2$, then $Y \subseteq H$.

\emph{Case 5.}  If $Y \notin \cL$ and $H \neq K^1, K^2$, then $Y$ has the form of $V$ in \eqref{E:reflat}.  
We may assume $H \cap W$ is a relative topoplane in $W$; it must be different from $H^1 \cap W$ and $H^2 \cap W$ since $\cH$ is simple.  
We work in the induced arrangement $\cH^W$.  In effect, that puts us in the situation where $W=X$, $Y=K^1$, and $Z = H^1 \cap H^2 = K^1 \cap K^2$.  Note that $Y \subseteq H^1 \cup H^2$.

Now there are several subcases depending on which of the intersections $H \cap H^i$ are void.  

\emph{Case 5a.}  If both are void, then $H\cap Y$ is empty.

\emph{Case 5b.}  Suppose one is void, say $H \cap H^1 \neq \eset = H \cap H^2$.  Then $H$, being disjoint from the relative topoplane $Z$ in $H^1$, lies in one half of $H^1$.  By choice of notation, $H \cap H^1 \subseteq H^1_+$.  

Now we make an argument that will show up again. $H \cap K^1 \subseteq K^1_+$, so $H \cap K^1 = H \cap H^1_+$, which (by Lemma \myref{L:nbdcross}) is a relative topoplane of $H^1_+$.  It follows that $H \cap K^1$ is a relative topoplane of $K^1_+$; we conclude that it is a relative topoplane of $K^1$.  This is what we needed to know in order to conclude that $\cH'$ is an arrangement of topoplanes.

\emph{Case 5c.}  Suppose that $H\cap H^1$ and $H \cap H^2$ are both nonempty.  Note that $H \not\supseteq Z$ by the simplicity of $\cH$.  Here we have two sub-subcases.

If $H \cap Z = \eset$, we can choose the notation so that $H \cap H^i \subseteq H^i_+$.  Then the argument of Case 5b implies that $H \cap K^i = H \cap H^i$, which is a relative topoplane both in $H$ and in $K^i$.

If $H \cap Z$ is not empty, then $V := H \cap Z$ is a relative topoplane in $Z$ and has codimension $3$.  $H \cap H^i$ has $V$ as a relative topoplane, so it is divided by $Z$ into $H \cap H^i_+$ and $H \cap H^i_-$, each of which is a relative topoplane in its half of $H^i$ and has as its boundary $H \cap Z$.  Now, 
$$
H \cap K^1 = (H \cap H^1_+) \cup (H \cap H^2_-) \cup (H \cap Z).
$$
In the right-hand side, the first part is a relative topoplane of $K^1_+$; the second part is a relative topoplane of $K^1_-$, and the last part is the boundary of each of the previous parts.  Thus, $H \cap K^1$ is a relative topoplane of $K^1$.  That is what we needed to show.

That ends the cases.  
To conclude the proof we observe that $\cH'$ has fewer noncrossing pairs of topoplanes than $\cH$, just as in Theorem \myref{T:planar}.  
By continuing with half-topoplane recombination we get an affine topoplane arrangement.
\end{proof}

\section{Topoplanes vs.\ pseudohyperplanes}

An \emph{arrangement of pseudospheres} in the $n$-sphere $S^n$ is a finite set $\cS$ of subspaces such that
\begin{itemize}
\item each $S \in \cS$ is a pseudosphere in $S^n$, i.e., $(S^n,S) \cong (S^n,S^{n-1})$ (where we think of $S^{n-1}$ as the equator of $S^n$) and $S$ is centrally symmetric in $S^n$,
\item the intersection of any subclass of $\cS$ is a topological sphere (which is necessarily again centrally symmetric), and
\item for any $\cS'\subseteq \cS$ and $S \in \cS\setm\cS'$, either $\bigcap\cS' \subseteq S$ or $S \cap \bigcap\cS'$ is a pseudosphere in $\bigcap\cS'$.
\end{itemize}
(It is known that every region is an open cell and its closure is a closed cell \cite{EM}.)
By identifying opposite points of $S^n$ we get a \emph{projective pseudohyperplane arrangement} $\cP$ in the real projective space $\bbP^n$.  
If we remove one pseudohyperplane $H_0 \in \cP$ from the arrangement and the space, and take the arrangement $\cA := \{ H\setm H_0 : H \in \cP, H \neq H_0 \}$ in $X := \bbP^n \setm H_0$, we have an \emph{affine pseudohyperplane arrangement}.  It is clearly an arrangement of topoplanes.  
We call a topoplane arrangement \emph{projectivizable} if it is homeomorphic to an arrangement constructed in this way, and more specifically we call it the \emph{affinization} of $\cP$.  
%
(See \cite[Chapter 5]{OM} for all facts about pseudosphere arrangements and \cite[Chapter 6]{OM} for projective pseudoline arrangements.)

There are several ways in which topoplane arrangements can be more complicated than affine pseudohyperplane arrangements.  In the analysis the concept of parallelism is important.  We define two topoplanes to be \emph{parallel} if they are disjoint.

\begin{lem}\mylabel{L:parallelequiv}
If a topoplane arrangement is projectivizable then it is affine and parallelism is an equivalence relation on topoplanes.
\end{lem}

\begin{proof}
It is easy to see from the known structure of pseudosphere, or projective pseudohyperplane, arrangements that $\cA$ is affine.  

Suppose $\cA$ is projectivizable.  Parallel topoplanes $H$ arise only from projective pseudohyperplanes $H_\bbP$ that meet at infinity.  If $H \parallel H' \parallel H''$, then $H_\bbP \cap H'_\bbP = Y$, a pseudohyperplane contained in the infinite hyperplane, and $H'_\bbP \cap H''_\bbP = Y$ also.  Thus, $H$ and $H''$ are parallel.
\end{proof}

\myexam{Disconnection}
The first way to get an unprojectivizable arrangement is by its being disconnected and not having all its topoplanes parallel.  We call a topoplane arrangement \emph{connected} if the union of its topoplanes (that is, the codimension-$1$ skeleton) is connected.  There are disconnected topoplane arrangements that are pseudohyperplane arrangements, indeed that are arrangements of true hyperplanes: take a finite family of parallel hyperplanes.  However, that is the only way.  It is just the opposite with topoplane arrangements.  Take any two topoplane arrangements $\cH_1$ and $\cH_2$ in two copies of $\bbR^n$.  In an unbounded region $R$ of $\cH_1$ find an open topological $n$-ball that extends to infinity.  By identifying this ball with $\bbR^n$ we can embed $\cH_2$ topologically inside $R$.
This gives a new topoplane arrangement $\cH := \cH_1 \cup \cH_2$ in $\bbR^n$ whose connected components are the components of $\cH_1$ and of $\cH_2$; in particular, assuming neither original arrangement was empty, the union is disconnected.  

\begin{prop}\mylabel{P:disconnected}
If $\cH_1$ has a pair of intersecting topoplanes, $\cH$ is not projectivizable.  
\end{prop}

\begin{proof}
The topoplanes in $\cH_1$ are parallel to those in $\cH_2$.  For $\cH$ to be projectivizable, parallelism must be an equivalence relation, so all the topoplanes are pairwise disjoint.  But this contradicts the assumption.
\end{proof}
\end{exam}

\myexam{The plane}\mylabel{X:planar}
In two dimensions nonequivalent parallelism is the only obstruction to being the affine part of a projective pseudoline arrangement.  (A \emph{pseudoline} is a pseudohyperplane in dimension $2$.)

\begin{thm}\mylabel{T:2dim}
An affine topoline arrangement in $\bbR^2$ is projectivizable if and only if parallelism in $\cA$ is an equivalence relation.
\dateadded{26 Aug 2006}
\end{thm}

\begin{proof}
The forward implication is obvious because topolines in the affinization are parallel if and only if they meet in a point at infinity.

For the converse, take a topoline arrangement $\cA$.  
Suppose it is affine and parallelism is an equivalence relation.  
Take a circle $C$ so large that all the intersection points as well as the other bounded faces of $\cA$ are inside $C$.  (If there is a topoline that is disjoint from all other topolines, imagine that it has a fictitious ``intersection point'' in the following discussion; that serves to make sure part of the topoline is inside $C$.)  
Each topoline $L^i \in \cA$ has two unbounded $1$-faces, which we arbitrarily label $L^i_+$ and $L^i_-$ and call the \emph{ends} of $L^i$.  Let $W^i_\epsilon$ be the first point on $L^i_\epsilon$, going from its finite end toward infinity, that lies on $C$.  We call the part of $L^i$ that extends from $W^i_\epsilon$ to infinity, away from the bounded part of $L^i$, the positive or negative \emph{tail} of $L^i$.  

To prove the theorem we replace the tails by new tails such that the positive tails of parallel topolines approach the same point at infinity, and the negative tails approach that point from the other side of infinity.  The rest of the proof explains a way to do that.

\begin{figure}[htdp]
\begin{center}
\includegraphics{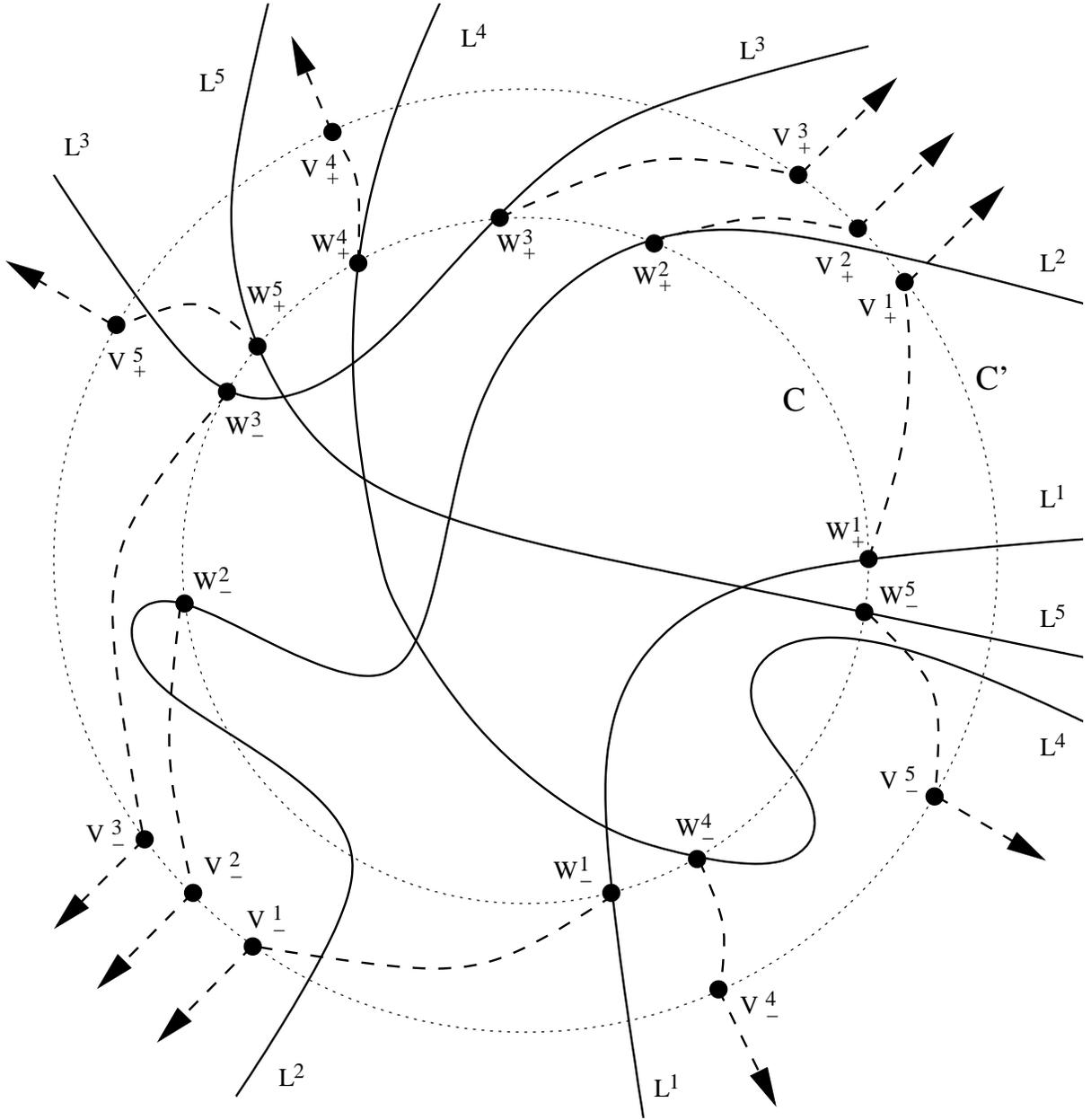}
\end{center}
\caption{The construction in the proof of Theorem \myref{T:2dim}, characterizing projectivizability of planar arrangements.}
\label{F:planar}
\end{figure}%

The points $W^i_\epsilon$ lie on $C$ in a cyclic order that is the same order in which the ends of the topolines appear outside $C$.  (The cyclic order of ends is well defined because there are no crossings outside $C$.)  
We show that the points of parallel topolines form two opposite consecutive groups.  
Suppose that $L^1 \parallel L^2$, and sign the $W$ points so their cyclic order is $W^1_+, W^2_+, W^2_-, W^1_-$.  Now suppose $W^3_+$ comes between $W^1_+$ and $W^2_+$.  If $L^3$ intersects $L^1$ it also intersects $L^2$, by transitivity of parallelism; but since $L^3_+$ is disjoint from $L^1$ and $L^2$, that forces the bounded faces in $L^3$ to intersect $L^1$ or $L^2$ twice, which is impossible.  Therefore, $L^3$ is parallel to $L^1$ and $L^2$ and, clearly, $W^3_-$ lies between $W^3_-$.  Thus, the $W$ points of a parallel class $L^1,\ldots,L^k$ appear in two consecutive groups along $C$, namely (in cyclic order around $C$) $W^1_+,\ldots,W^k_+,S_+,W^k_-,\ldots,W^1_-S_-$, where $S_\epsilon$ is the set of $W^i_\epsilon$ points of all other topolines $L^i$, since each of those $L^i$ crosses all of $L^1,\ldots,L^k$.  Let us call the points $W^i_\epsilon$ of each group, but with fixed $\epsilon$, \emph{equivalent points}.  Changing the signs of the points in an equivalence class gives the \emph{opposite} class.

Choose a larger circle $C'$ concentric with $C$ and points $V^i_\epsilon$ on $C'$ in the same cyclic order as the $W^i_\epsilon$, and give them the same equivalence relation.  Pick the $V$ points so that those in one equivalence class are close together.  Furthermore, if $V_+$ and $V_-$ denote the midpoints of the arcs containing an equivalence class and its negative, the points should be chosen so $V_+$ and $V_-$ are diametrically opposed.  Draw nonintersecting curves in the annulus bounded by $C$ and $C'$ that connect corresponding $W$ and $V$ points.

For each equivalence class of $V$ points, choose the direction $d$ that extends from its midpoint $V_\epsilon$ radially away from the center of $C'$.  Draw rays from each point in the equivalence class in the direction $d$.  Now we replace each topoline $L^i$ by the curve made up of the part of $L^i$ that is not in the tails, together with the two curves from $W^i_\epsilon$ to $V^i_\epsilon$ and the rays emanating from the two points $V^i_\epsilon$.  
By the rule for choosing midpoints, opposite classes have opposite directions.  Since the points of each class are close together, the rays are entirely outside $C'$ and therefore do not intersect each other or any of the curves from $W$ points to $V$ points or any of the parts of the original topolines other than their tails.  Thus, the new topolines form an arrangement $\cA'$ that has the same intersection points (and all bounded faces) as the original ones.  
It is clear that $\cA'$ is homeomorphic (indeed isotopic) to $\cA$.

Moreover, the topolines of $\cA'$ have the property that parallels approach the same point at infinity while nonparallels do not.  Furthermore, the opposite ends of the new topolines approach the same point at infinity, but from opposite directions.  
Thus, we can add the infinite line to get a projective pseudoline arrangement $\cP$ from which $\cA'$ is derived by affinization; and $\cA$ is homeomorphic to $\cA$, so $\cA$ is projectivizable.
\end{proof}
\end{exam}

\myexam{Connected, affine, but not projectivizable}\mylabel{X:noproj2}
To get a simple example of an affine topoline arrangement that is not projectivizable, take the four topolines  $x_1 = -1$, $x_1 = 1$, $x_2 = 1$, and the bent line $\{ x: x_1x_2=0 \text{ and } x_1, x_2 \geq 0 \}.$  In this example parallelism is obviously not transitive.  One can even omit the horizontal line, but it  is what makes the arrangement connected.
\end{exam}

In higher dimensions, which affine topoplane arrangements are projectivizable remains mysterious.  Is intransitivity of parallelism the only obstruction?

\section{No weaker definition}

Examples show that our definition of an arrangement of topoplanes cannot be simplified in some tempting ways.  The essential property of flats for the proof of Equation \eqref{E:count} is that a flat $Y$ has a rank, $r(Y)$, in the intersection semilattice and its Euler characteristic is $(-1)^{r(Y)}$.  The natural way to ensure this is to require that $Y$ have codimension equal to its rank, and be homeomorphic to $\bbR^{\dim Y}$.  The essential property of regions is that each open region be a cell; this seems to require that a flat be a topoplane in each flat that it covers.  However, that alone is not enough; and this is not the only natural idea for simplifying the definition that does not work.

\myexam{Pair intersection}\mylabel{X:pairmeet}
For instance, it would be much simpler if it were sufficient that pairs topoplanes intersect in a relative topoplane of each.  Here is a counterexample consisting of three topoplanes, each pair intersecting in a relative topoplane, but the intersection of all three being neither a relative topoplane nor of the correct dimension.  
In $\bbR^3$ let $H_1$ be the plane $x_1=-x_2$ and let $H_2$ be the plane $x_1=x_2$.  For $H_3$ we use the surface defined by 
$$
x_2 = \begin{cases} x_3-1 &\text{if } x_3\geq 1, \\ 0 &\text{if } x_3\in[-1,1], \\ x_3+1 &\text{if } x_3\leq-1.\end{cases}
$$
Each $H_i \cap H_j$ is a straight line or a broken line that divides $H_i$ and $H_j$ into two parts, but the intersection of all three topoplanes is the line segment $\{ (0,0,x_3): -1\leq x_3 \leq 1 \}$.
\end{exam}

\myexam{Flat intersection}\mylabel{X:badmeet}
One might still hope it would be sufficient that, if a flat $Y$ covers a flat $Z$, then $Z$ is a relative topoplane of $Y$.  
(In $\cL$ we say $Y$ \emph{covers} $Z$ if $Y > Z$---that is, $Y \subset Z$---and there is no other element in between them.)  
Another example of three topoplanes shows that this is too weak to give us an arrangement of topoplanes.  In $X = \bbR^3$ take the two halves of the cone $x_2^2+x_3^2=1$, one opening to the right and the other to the left, to be $H_1$ and $H_2$.  Let $H_3$ be a plane tangent to the cone in a line $W$ and let $Z := $ the origin.  Setting $\cH := \{H_1,H_2,H_3\}$, the intersection poset is $\cL = \{\bbR^3,H_1,H_2,H_3,W,Z\}$.  This satisfies the covering property but it is not a topoplane arrangement because $H_1 \cap H_2$ is not a topoplane in $H_1$.
\end{exam}

Still, none of these counterexamples applies to arrangements of topolines; for them, it is sufficient to require only that the intersection of any two topolines be void or a point.  It is also sufficient to require that for any covering pair $Y,Z$, $Z$ is a relative topoline in $Y$, except that one must forbid the case of a single flat that is a point.  (These facts are obvious.)



\begin{thebibliography}{9}

\bibitem{OM} A.\ Bj{\"o}rner, M.\ Las Vergnas, B.\ Sturmfels, N.\ White, and G.M.\ Ziegler, 
\emph{Oriented Matroids}.
Encyclopedia Math.\ Appl., Vol.\ 46. 
Cambridge Univ.\ Press, Cambridge, 1993.
MR 95e:52023.  Zbl.\ 773.52001.
Second ed., 1999.
MR 2000j:52016.  Zbl.\ 944.52006.

\bibitem{EM} Jack Edmonds and Arnaldo Mandel, 
Topology of oriented matroids.
Abstract 758-05-9.
\emph{Notices Amer.\ Math.\ Soc.} {\bf 25} (1978), A-510.
Abstract of \cite{Mandel}.

\bibitem{FL} Jon Folkman and Jim Lawrence, 
Oriented matroids. 
\emph{J.\ Combin.\ Theory Ser.\ B} \textbf{25} (1978), 199--236.
MR 81g:05045.  Zbl.\ 325.05019.

\bibitem{Mandel} Arnaldo Mandel, 
\emph{Topology of Oriented Matroids}.
Doctoral dissertation with and supervised by Jack Edmonds, 
University of Waterloo, 1982.

\bibitem{WW} Michelle L.\ Wachs and James W.\ Walker, 
On geometric semilattices. 
\emph{Order} {\bf 2} (1986), 367--385.
MR 87f:06004.  Zbl.\ 589.06005.

\bibitem{CATD} Thomas Zaslavsky,
A combinatorial analysis of topological dissections.
\emph{Advances in Math.} \textbf{25} (1977), 267--285. 
MR 56 \#5310.  Zbl.\ 406.05004.

\end{thebibliography}
\end{document}